\documentclass[10pt]{article}
\usepackage{amssymb}
\usepackage{amsmath}

\usepackage{amssymb}
\usepackage{amsfonts}
\usepackage{amscd}
\newtheorem{theorem}{Theorem}

\newtheorem{definition}[theorem]{Definition}

\newtheorem{lemma}[theorem]{Lemma}

\newtheorem{proposition}[theorem]{Proposition}

\newenvironment{proof}[1][Proof]{\textbf{#1.} }{\ \rule{0.5em}{0.5em}}

\newcommand\R{\mathbb{R}}
\newcommand\g{\mathfrak{g}}

\newcommand\C{\mathbb{C}}
\newcommand\Z{\mathbb{Z}}

\newcommand\h{\mathfrak{h}}
\newcommand\m{\mathfrak{m}}
\newcommand\gm{\Gamma}
\newcommand\Aut{\mathrm{Aut}}
\newcommand\wg{\widehat{\Gamma}}
\newcommand\wP{\widehat{P}}
\newcommand\z{\mathbb{Z}_2\times \mathbb{Z}_2}

\newcommand\soo{\mathrm{so}}
\newcommand\spp{\mathrm{sp}}
\newcommand{\cj}[2]{#1^{-1}\!#2#1}
\newcommand{\su}[1]{\mathrm{Supp}\,#1}

\newcommand{\sll}{\mathrm{sl}}

\newcommand{\ad}{\mathrm{ad}\,}

\newcommand{\gl}[1]{\mathrm{gl}(#1)}
\newcommand{\gll}{\mathrm{gl}}
\newcommand{\tr}[1]{#1^t}
\newcommand{\as}[1]{#1^{\ast}}
\newcommand{\inv}[1]{#1^{-1}}

\newcommand{\ve}{\varepsilon}
\newcommand{\vp}{\varphi}

\begin{document}

\title{$\z$-symmetric spaces}
\author{Yuri Bahturin and Michel Goze}
\maketitle
\begin{abstract}
The notion of a $\Gamma $-symmetric space is a generalization 
of the  classical notion of a symmetric space, where a general
finite abelian group $\Gamma $ replaces the group $\Z_2$. The case $\Gamma =\Z_k$ has also been studied, from the algebraic point of view by V.Kac \cite{VK} and from the point of view of the differential geometry by Ledger, Obata \cite {LO}, Kowalski \cite{K} or Wolf - Gray \cite{WG} in terms of
$k$-symmetric spaces. In this case, a $k$-manifold is an homogeneous reductive space and the classification of these varieties is given by the corresponding classification of graded Lie algebras. The general notion of a $\Gamma $-symmetric space was introduced by R.Lutz in \cite{L}. We approach the classification of such spaces in the case $\gm=\z$ using recent results (see \cite{BShZ}) on the classification of complex $\z$-graded simple Lie algebras.
\end{abstract}
\medskip


\section{Introduction}\label{sI}

 
A \emph{symmetric space} is a homogeneous space $M=G/H$ where $G$ is a connected Lie group with an involutive automorphism
$\sigma $ and $H$ a closed subgroup which lies between the subgroup of all fixed points of $\sigma $ 
and its identity component.
This automorphism $\sigma $  induces an involutive {\it diffeomorphism} $\sigma _0$ of $M$ such that 
$ \sigma _0(\pi (x))=\pi (\sigma (x))$ for every $x \in G$ where $\pi :G\rightarrow G/H$ is the canonical projection. It 
also induces an automorphism $\gamma $ on the Lie algebra $\g$ of $G$. This automorphism satisfies $\gamma ^2=Id$, hence
is diagonalizable and the Lie algebra $\g$ of $G$ admits a $\Z_2$-grading, $\g=\g_0 \oplus \g_1$ where $\g_0$ and 
$\g_1$ are the 
eigenspaces of $\sigma $ corresponding to the eigenvalues $1$ and $-1$. Conversely, every $\Z_2$-grading
$\g=\g_0 \oplus \g_1$ on a Lie algebra $\g$
makes it into a \emph{symmetric Lie algebra}, that is a triple $(\g,\g_0,\g  )$ where $\gamma $ 
is an involutive automorphism of $\g$ such that
$\gamma (X)=X$ if and only if $X \in \g_0$ and $\gamma (X)=-X$ for all $X \in \g_1$. If $G$ is a connected simply 
connected Lie group with Lie algebra $\g$, then
$\gamma $ induces an automorphism $\sigma$ of $G$ and for any subgroup $H$ lying between 
$G^{\sigma }=\{x \in G, \sigma (x)=x\}$
and the identity component of $G^{\sigma }$, $(G/H,\sigma )$ is a symmetric space. In the Riemannian case, 
$H$ is compact
and $\g$ admits an orthogonal symmetric decomposition, that is the Lie group of linear transformations of $\g$
generated by $\ad_{\g}H$ is compact. As a result, the study is reduced to the effective irreducible case and $\g$ is 
semi-simple. 

In this paper we will look at more general $\Gamma $-symmetric homogeneous spaces. They were first introduced by R. Lutz in \cite{L}, and A. Tsagas at a workshop in Bucharest. We propose here to develop the corresponding algebraic structures
and to give, using the results on complex simple Lie algebras graded by finite abelian groups \cite{BShZ},  \cite{BZA},\cite{BZ}, (see also \cite{P1},\cite{P2}), the algebraic classification of $\z$-symmetric spaces $G/H$ whose associated Lie algebra $\g$ is classical simple. 


\section{Group graded Lie algebras}\label{sggla}


\subsection{Definition}
\begin{definition}
Let $P$ be a group with identity element $1$. A  Lie algebra $\g$ over a field $F$ is graded by  $P$ if
$$\g=\bigoplus _ {p\in P } \ \g_p$$ with $$[\g_{p},\g_ {q}]\subset \g _{pq}$$ for all $p,q\in P$.
\end{definition}

\begin{definition}\label{deq}
Given two $P$-gradings 
$$
\g=\bigoplus_{p\in P}\g_{p}\mbox{ and 
}\g=\bigoplus_{p\in P}\g_{p}^{\prime}
$$
of an algebra $\g$ by a group $P$ 
we call them \textit{equivalent} if there exists an automorphism $\theta$ of $\g$ such that 
$\g_{p}^{\prime}=\theta(\g_{p})$, for all $p\in P$.
\end{definition}

An important subset of the grading group is defined by the following.

\begin{definition}
Given a grading as above, the set 
$$\su\g=\{p\in P\, |\,\g_p\neq\{ 0\}\}
$$
is called the \emph{support} of the grading.
\end{definition}
It has been established in \cite{P0} (see also \cite{BZ}) that if $\g$ is complex simple then any two elements in the support of the grading commute. So one can always restrict oneself to the case of abelian groups. In this paper we restrict ourselves to finite abelian grading groups $P$.


\subsection{Action of the dual group}\label{s12}


Let $\Gamma =\wP$ be the dual group associated to $P$, that is, the group of characters $$\gamma :P \longrightarrow \C^\star $$ of $
P $. If we assume that  a Lie algebra $\g$ is $P$-graded then we obtain a natural action of  $\gm$  by linear transformations on $\g\otimes\C$ if for any homogeneous elements
 $X\in \g_p$ we set $\gamma(X)=\gamma(p)X$. Since for $X\in\g_p$ and $Y\in\g_q$ we have $[X,Y]\in\g_{pq}$, it follows that
 \begin{equation}\label{ega} \gamma([X,Y]=\gamma(pq)[X,Y]=[\gamma(p)X,\gamma(q)Y]=[\gamma(X),\gamma(Y)],
 \end{equation}
that is, $\Gamma$ acts by Lie automorphisms on $\g$. In this case there is a canonical homomorphism  
\begin{equation}\label{ech}
\alpha:\gm\rightarrow \mathrm{Aut}(\g\otimes\C)\mbox{ given by } \alpha(\gamma)(X)=\gamma(X). 
\end{equation}
If for any $p\in\su{\g}$ we have $p^2=1$ then the action is defined even on $\g$ itself and the above homomorphism maps $\gm$ onto a subgroup of $\mathrm{Aut}\,\g$.

Conversely, suppose there is a homomorphism $\alpha:\gm\rightarrow \mathrm{Aut}\,\g$, for a finite abelian group $\gm$. Then $\gm$ acts on $\g$, hence on $\g\otimes\C$ by automorphisms if one sets $\gamma(X)=\alpha(\gamma)(X)$. This action extends to $\g\otimes\C$ and yields a $P$-grading, $\gm=\wP$, of $\g\otimes\C$ by subspaces $(\g\otimes\C)_p$, for each $p\in P$, defined as follows: $$(\g\otimes\C)_p=\{X\,|\,\gamma(X)=\gamma(p)X\}.$$
That $[\g_p,\g_q]\subset\g_{pq}$ easily follows from (\ref{ega}). Now the vector space decomposition $\g=\bigoplus_{p\in P}\g_p$ is just a standard weight decomposition under the action of an abelian semisimple group of linear transformations over an algebraically closed field. Again, if we have $\alpha(\gamma)^2=1$ for any $\gamma\in\gm$, then we have a $P$-grading on $\g$ itself.

Explicitly, let us assume that $\g$ is a complex Lie algebra and $K$ a finite abelian subgroup of $\Aut(\g)$. One can write $K$ as $K=K_1\times ...\times K_p$
where $K_i$ is a cyclic group of order $r_i$. Let $\kappa _i$ be a generator of $K_i$. The automorphisms $\kappa _i$
satisfy:
$$
\left\{
\begin{array}{l}
\kappa _i ^{r_i}=Id,\\
\kappa _i\circ \kappa _j=\kappa _j\circ \kappa _i,
\end{array}
\right.
$$
for all $i,j=1,...,p.$ These automorphisms are simultaneously diagonalizable. If $\xi _i$ is a primitive root 
of order $r_i$ of the unity, then the eigenspaces
$$\g_{s_1,...,s_p}=\{X \in \g \ {\mbox{\rm such \ that}} \ \kappa _i(X)=\xi _i^{s_i}X, \ i=1,...,p\}$$
give the following grading of $\g$ by $P=\Z_{r_1}\times ...\times \Z_{r_p}$ :
$$\g=\bigoplus _{(s_1,...,s_p) \in \Z_{r_1}\times ...\times \Z_{r_p}} \ \g_{s_1,...,s_p}.$$
We can summarize some of what was said above as follows.

\begin{proposition}\label{p4}
Let $P$ be a finite abelian group and $\gm=\wP$, the group of complex characters of $P$. 
\begin{enumerate}
\item[{\rm(a)}]
A complex Lie algebra $\g$ is $P$-graded if and only if the dual group
$\Gamma $ maps homomorphically onto a finite abelian subgroup of $\Aut(\g)$, by a canonical homomorphism $\alpha$ described in (\ref{ech}). 
\item[(b)] A real Lie algebra $\g$ is $P$-graded, with $p^2=1$ for each $p\in\su{\g}$, if and only if there is a homomorphism $\alpha:\gm\rightarrow\Aut(\g)$ such that $\alpha(\gamma)^2=\mathrm{id}_{\g}$ for any $\gamma\in\gm$.
\item[(c)] In both cases above, $\su{\g}$ generates $P$ if and only if the canonical mapping $\alpha$ has trivial kernel, that is, $\gm$ is isomorphic to a (finite abelian) subgroup of $\Aut(\g)$.
\end{enumerate}
\end{proposition}

\begin{proof} We need only to comment on (c). If $\Lambda\subset\gm$ then by $\Lambda^\perp$ we denote the set of all $p\in P$ such that $\lambda(p)=1$ for all $\lambda\in\Lambda$. Similarly we define $Q^\perp$ for any $Q\subset P$. We have that $\Lambda^\perp$ and $Q^\perp$ are always subgroups in $P$ and $\gm$, respectively. If $\Lambda$ and $Q$ are subgroups then $|\Lambda|\cdot|\Lambda^\perp|=|Q|\cdot|Q^\perp|=|\gm|=|P|$. We claim that if $\Lambda=\mathrm{Ker}\,\alpha$ then the subgroup generated by $\su{\g}$ is $P=\Lambda^\perp$. This follows because for any $p\notin Q$ there is $\lambda\in\Lambda$ such that $\lambda(p)\neq 1$. If $\g_p\neq \{ 0\}$ and $0\neq X\in\g_p$ then $\lambda(X)=\lambda(p)X\neq X$ and $\lambda\notin\mathrm{Ker}\,\alpha$. Conversely, let $\su{\g}\subset T$, where $T$ is a proper subgroup of $Q= \mathrm{Ker}\,\alpha^\perp$. Then $T^\perp$ properly contains $\Lambda$ and for any $\gamma\in T^\perp\setminus\Lambda$ and any $X\in \g_t$, $t\in T$, we have $\gamma(X)=\gamma(t)X=X$. Since all such $X$ span $\g$, we have that $\gamma\in\mathrm{Ker}\,\alpha=\Lambda$, a contradiction. Thus $\su{\g}$ must generate the whole of $Q$, proving (c).
\end{proof}

\medskip

 
 \subsection{Examples}\label{ssE}
 

 1. The gradings of classical simple complex Lie algebras by finite abelian groups have been described in \cite{BSZ}, \cite{BShZ}, \cite{BZ}, and \cite{DV} (see also \cite{P1}, \cite{P2}). We will use this classification in the case $\gm=\z$ to obtain some classification - type results in the theory of $\z$-symmetric spaces.

\medskip

 2. In the non simple case the study of gradings is more complicated. Consider, for example, the nilpotent case. In
contrast to the simple case, there is no classification of these Lie algebras, except in the dimensions up to $7$, \cite{GR1}. Even then one has to distinguish between two classes of complex nilpotent Lie algebras:

\begin{enumerate}

\item[(a)] The non characteristically nilpotent Lie algebras. These Lie algebras admit a non trivial abelian subalgebra
of  the Lie algebra $\mathrm{Der}(\g)$ of derivations such that the elements are semisimple. In this case $\g$ is graded by the
roots.
 
\item[(b)] The characteristically nilpotent Lie algebras. Every derivation is nilpotent and we do not have root decompositions. 
Nevertheless, these nilpotent Lie algebras can be graded by groups. For example, the following nilpotent Lie algebra,
denoted by $(\mathfrak{n}_7^3)$ in \cite{GR1} and given by
$$
\left\{
\begin{array}{l}
\lbrack X_1,X_i \rbrack =X_{i+1} , \ i=2,...,6\\
\lbrack X_2, X_3 \rbrack = X_5+X_7 \\
\lbrack X_2, X_4 \rbrack = X_6 \\
\lbrack X_2, X_5 \rbrack = X_7 \\
\end{array}
\right.
$$
is characteristically nilpotent and admits the following grading
$$\mathfrak{n}_7^3=\g_0 \oplus \g_1$$
where $\g_0$ is the nilpotent subalgebra generated by $\{X_2,X_4,X_6\}$ and $\g_1$ is the $\g_0$-module
generated by $\{X_1,X_3,X_5,X_7\}$. On the other hand, the nilpotent Lie algebras 
$\mathfrak{n}_7^4$ and $\mathfrak{n}_7^5$ do not admit nontrivial group gradings. 

\end{enumerate}

 
\section{$\Gamma$ -symmetric spaces}\label{sgss}


\subsection{Definition}\label{ssD}


\begin{definition}
Let $\Gamma $ be a finite abelian group.  A homogeneous space $M=G/H$ is called \emph{$\Gamma$-symmetric } if
\begin{enumerate}
\item The Lie group $G$ is connected

\item The group $G$ is effective on $G/H$ (i.e the Lie algebra $\h$ of $H$ does not contain a nonzero proper ideal of 
the Lie algebra $\g$ of $G$)

\item There is an injective homomorphism
$$ \rho: \gm \rightarrow \mathrm{Aut}\, G$$
such that if $G^{\gm }$ is the closed subgroup of all elements of $G$ fixed by $\rho (\gm)$
and $(G^{\gm})_e$ the identity component of $G^\Gamma $ then
$$(G^{\gm})_e\subset H \subset G^{\gm}.$$
\end{enumerate}
\end{definition}

Obviously, in the case of $\gm=\mathbb{Z}_2$ we obtain ordinary symmetric spaces.

\medskip

We denote by $\rho _\gamma $ the automorphism $\rho (\gamma )$ for any $  \gamma  \in \gm$.  If $H$ is connected, we have
 $$
\left\{
\begin{array}{l}
\rho _{\gamma  _1}\circ \rho _{\gamma _2}=\rho _{\gamma _1 \gamma _2}, \\
\\
 \rho _{\ve}=Id  \\
\\
\rho _\gamma(g)=g \ ,  \forall \gamma\in \gm  \Longleftrightarrow g\in H.
\end{array}
\right.
$$
where $\ve$ is the identity element of $\gm $. If $\Gamma =\Z_2$ then a $\Z_2$-symmetric space is a symmetric space,
if $\Gamma =\Z_p$ we find again the $p$-manifolds in the sense of Ledger-Obata \cite{LO}.
 

\subsection{$\wg$-grading of the Lie algebra of $G$}\label{sswlg}


Let $M=G/H$ be a $\Gamma $-symmetric space. Each automorphism $\rho _\gamma$ of $G$, $\gamma\in\gm$,  induces an automorphism of $\g$,
denoted by $\tau_\gamma$ and given by $\tau_\gamma =(T\rho _\gamma )_e$ where $(Tf)_x$ is the tangent map
of $f$ at the point $x$.

\begin{lemma}\label{l6}
The map $\tau :\gm \longrightarrow \Aut(\g)$ given by 
$$\tau (\gamma)=(T\rho _\gamma)_e$$
is an injective homomorphism.
\end{lemma}

\begin{proof} Let $\gamma _1,\gamma _2$ be in $\Gamma $. Then $ \rho _{\gamma _1}\circ \rho _{\gamma _2}=
\rho _{\gamma _1\gamma _2}$. It follows that $ (T \rho _{\gamma _1})_e\circ(T \rho _{\gamma _2})_e=
(T \rho _{\gamma _1}\rho _{\gamma _1})_e=(T\rho(\gamma_1\gamma_2))_e$, that is, $\tau (\gamma_1\gamma _2 )=\tau (\gamma _1)\tau (\gamma _2)$. Now
let us assume that $\tau (\gamma )=Id_{\g}$. Then $(T\rho _\gamma )_e=Id=(T\rho _{\ve} )_{e}$. But
$\rho _\gamma $ is uniquely determined by the corresponding tangent automorphism of $\g$. Then $\rho _\gamma =\rho _\ve$
and $\gamma =\ve$.
\end{proof}

\medskip

From this lemma we derive the following.

\begin{proposition}\label{p7}
If $M=G/H$ is a $\Gamma $-symmetric space then the complex Lie algebra $\g_{\C}=\g\otimes \C$ where $\g$ is the
Lie algebra of $G$ is $\wg $-graded and if $\Gamma =\Z_2^k$ then the real Lie algebra $\g$ of $G$ 
is $\wg $-graded. The subgroup of $\wg$ generated by the support of the grading is $\wg$ itself.
\end{proposition} 
\begin{proof} Indeed, by Lemma \ref{l6}, $\alpha: \Gamma \rightarrow \Aut\,\g$ is an injective homomorphism, so all our claims follow by Proposition \ref{p4}.
\end{proof}

For convenience, recall that if a finite abelian group $P$ is such that $\wg=P$ and $\alpha$ is a canonical homomorphism introduced in (\ref{ech}) then the components of the grading are given by the following equation.
\begin{equation}\label{egc}
\g_p =\{X \in \g, \,|\, \alpha (\gamma)(X)=\gamma(p)X, \ \forall  \gamma \in  \Gamma \}.
\end{equation}


\subsection{$\Gamma$-symmetric spaces and graded Lie algebras}\label{ssgss}


To study $\Gamma$-symmetric spaces, we need to start with the study  of $P$-graded Lie algebras where $\gm=\wP$. But in a general case
if $G$ is a connected Lie group corresponding to $\g$, the $P$-grading of $\g$ or $\g_{\C}$ does not necessarily give a 
$\Gamma $-symmetric space $G/H$. Some examples are given in \cite{Be}, even in  the symmetric case. Still,
if $G$ is simply connected, $\Aut(G)$ is a Lie group isomorphic to $\Aut(\g)$.
 
\begin{proposition}\label{p8}
Let $P=\Z_2^k$, with the identity element $\varepsilon$, $\wg=P$, and $\g$ a real $P$-graded Lie algebra such that the subgroup generated by $\su{\g}$ equals $P$ and the identity component $\h=\g_{\epsilon}$ of the grading does not contain a nonzero ideal of $\g$. If $G$ is a connected simply connected Lie group with the Lie algebra $\g$ and $H$ a Lie subgroup associated with $\h$, then the homogeneous space $M=G/H$ is a $\gm$-symmetric space.
\end{proposition} 

\begin{proof} By Proposition \ref{p4}, there is an injective homomorphism $\alpha:\gm\rightarrow\Aut\,\g$ defined by this grading. The subgroup $\alpha(\gm)$ of $\Aut\,\g$ is isomorphic to $\gm$. Choosing, for each $\alpha(\gamma)$, a unique automorphism $\rho(\gamma)$ of $G$ such that $(T(\rho(\gamma))_e=\tau(\gamma)$ we obtain an injective homomorphism  $\rho:\gm\rightarrow\Aut\, G$, making $G/H$ into a $\gm$-symmetric space.
\end{proof}

\medskip

Motivated by Propositions \ref{p4}, \ref{p7} and \ref{p8}, we introduce the following.

\begin{definition}\label{dls}
Given a real or complex Lie algebra $\g$, a subalgebra $\h$ of $\g$ and a finite abelian subgroup $\gm\subset\Aut\,\g$, we say that $(\g/\h,\gm)$ is a \emph{local $\gm$-symmetric space} if $\h=\g^{\gm}$, the set of all fixed points of $\g$ under the action of $\gm$. We call $\h$ the \emph{isotropy subalgebra} of $(\g/\h,\gm)$.
\end{definition}

Any $\Gamma$-symmetric space gives rise to a local $\gm$-symmetric space and, in the case of connected simply connected groups, the converse is also true. If $\gm=\mathbb{Z}_2^k$ then $(\g/\h,\gm)$ is a local $\gm$-symmetric space if and only if $\g$ is $P$-graded, where $P=\wg$, and the isotropy subalgebra $\h$ is the identity component of the grading. If $\gm$ is a more general finite abelian group then the grading by $P=\wg$ arises only on the complexification $\g\otimes\C$ and still $\h\otimes\C$ is both the set of fixed points of $\gm$ and the identity component of the grading. Again, the study of local complex $\gm$-symmetric spaces amounts to the study of $P$-graded Lie algebras, where $P=\wg$.


\subsection{$\Gamma $-symmetries on the homogeneous space $M=G/H$}\label{ssgshs}


Given a $\Gamma $-symmetric space $(G/H,\Gamma)$ it is easy to construct, for each point $x$ of the
homogeneous space $M=G/H$, a subgroup
of the group $\mathrm{Diff}(M)$ of diffeomorphisms of $M$, isomorphic to $\Gamma $, which has  $x$ as an isolated fixed point. We denote by 
$\bar g$ the class
of $g \in G$ in $M$. If $e$ is the identity of $G$, $\gamma\in\gm$,
we set
$$s  _{(\gamma ,\bar e)}(\bar g)=\overline{\rho _\gamma (g)}.$$
If $\bar g$ satisfies $ s  _{(\gamma,\bar e)}(\bar g)=\bar g$ then $\overline{\rho_\gamma (g)}=\bar g$
that is $\rho _\gamma  (g)=gh_\gamma $ for $h_\gamma \in H$. Thus $h_\gamma =g^{-1}\rho _\gamma  (g)$. 
But $\Gamma\cong\wg $ is a finite abelian group. If $p_\gamma $ is the order of $\gamma $ then 
$\rho _{\gamma ^{p_\gamma }}=Id.$ Then 
$$h_\gamma ^{2}=g^{-1}\rho _\gamma(g)\rho _\gamma (g^{-1}) \rho _{\gamma ^2}(g)=g^{-1}\rho _{\gamma ^2}(g).$$
Applying induction, and considering $h^m \in H$ for any $m$, we have
$$ (h_\gamma )^m=g ^{-1}\rho _{\gamma ^m}(g).$$
For $m=p_\gamma $ we obtain
$$(h_\gamma)^{p_ \gamma}=e.$$ 
If $g$ is near the identity element of $G$, then $h_\gamma $ is also close to the identity and $h_\gamma ^{p_\gamma }=e$
implies $h_\gamma =e.$ Then $\rho _\gamma (g)=g$. This is true for all $\gamma \in \gm$ and thus $g \in H$. It follows that
$\overline g=\overline e$ and that the only fixed point
of $s  _{(\gamma ,\bar e)}$ is $\bar e$. In conclusion, the family 
$\{s  _{(\gamma ,\bar e)}\}_{\gamma \in \gm}$ of diffeomorphisms of $M$ satisfy
$$
\left\{
\begin{array}{l}
s  _{(\gamma _1 ,\bar e)}\circ s  _{(\gamma _2,\bar e)}=s  _{(\gamma _1\gamma _2 ,\bar e)}\\
s  _{(\gamma ,\bar e)}(\bar g)=\bar g, \forall \gamma \in \gm \Rightarrow \bar g=\bar e.
\end{array}
\right.
$$
Thus,
$$\gm  _{\bar e }=\{ s_{(\gamma ,\bar e)}, \ \gamma \in \gm \}$$ 
is a finite abelian subgroup of $\mathrm{Diff}(M)$ isomorphic to $\Gamma $, for which $\bar e$ is an isolated fixed point.

\medskip

 In another point $\bar {g_0} $ of $M$ we put
$$s  _{(\gamma ,\bar {g_0})}(\bar g)=g_0(s  _{(\gamma ,\bar e)})(g_0^{-1}\bar g).$$ As above, we can see that
$$
\left\{
\begin{array}{l}
s  _{(\gamma _1 ,\bar {g_0})}\circ s  _{(\gamma _2,\bar {g_0})}=s  _{(\gamma _1\gamma _2 ,\bar {g_0})}\\
s  _{(\gamma ,\bar {g_0})}(\bar g)=\bar g, \forall \gamma \in \gm \ \Rightarrow \bar g=\bar {g_0}.
\end{array}
\right.
$$
and
$$\Gamma  _{\bar {g_0} }=\{ s_{(\gamma ,\bar {g_0})}, \ \gamma \in \gm \}$$ 
is a finite abelian subgroup of $\mathrm{Diff}(M)$ isomorphic to $\Gamma $, for which $\bar g_0$ is an isolated fixed point.

Thus for each $\bar g\in M$ we have
a finite abelian subgroup $\Gamma _{\bar g}$ of $\mathrm{Diff}(M)$ isomorphic to $\Gamma $,  for which $\bar g$ is an isolated fixed point. 
 
\begin{definition}
Let $(G/H,\Gamma )$ be a $\Gamma $-symmetric space. For any point $x \in M=G/H$ the subgroup $\Gamma _x\subset\mathrm{Diff}(M)$
is called \emph{the group of symmetries} of $M$ at $x$. 
\end{definition}
 
\medskip

Since for every $x \in M$ and $\gamma \in \gm $, the map $ s _ {(\gamma ,x)} $ is a diffeomorphism of $M$
such that $s _ {(\gamma ,x)}(x)=x$, the tangent linear map
 $(Ts _ {(\gamma ,x)})_x$ is in $\mathrm{GL}(T_xM)$. 
For every $x \in M$, we obtain a linear representation
$$S _x:\gm \longrightarrow \mathrm{GL}(T_xM)$$
defined by
$$S_x(\gamma )=(T s _ {(\gamma ,x)})_x.$$ 
Thus for every $\gamma  \in \gm $ the map
$$S(\gamma ): M \longrightarrow T(M)$$
defined by $S(\gamma )(x)=S_x(\gamma )$ 
is a $(1,1)$-tensor on $M$ which satisfies:
\begin{enumerate}
\item the map $S(\gamma )$ is of class ${\mathcal{C}}^\infty $,

\item for every $x \in M$, 
$$\{X_x \in T_x(M)\,|\, S_x(\gamma )(X_x)=X_x, \forall \gamma \in \gm \}=\{0\}.$$
In fact, this last remark is a consequence of the property : 
$s _{(\gamma ,x)}(y)=y$ for every $\gamma $ implies $y=x$.
\end{enumerate}
\begin{definition}
Let $M$ be a real differential manifold and $\Gamma $ a finite abelian group. We note by $T_xM$ the tangent 
space to $M$ at the point $x$.

A $\Gamma $-symmetric structure on $M$ is given for all $x \in M$ by a linear representation of $\gm $ on the
vector space $T_xM$ 
$$\rho_x :\gm \longrightarrow \mathrm{GL}(T_xM)$$
such that 
\begin{enumerate}
\item For every $\gamma \in \gm $, the map $x \in M \longrightarrow \rho_x (\gamma )$ is of
class $\cal{C}^{\infty }$

\item For every $x \in M$, $\{X_x \in T_x(M)\,|\, \rho _{x}(\gamma )(X_x)=X_x, \ \forall \gamma \}=\{0\}$.
\end{enumerate}
\end{definition}
\begin{proposition}
If $(G/H,\Gamma)$ is a $\Gamma $-symmetric space, the family $\{S _x\}_{x \in M}$ is a $\Gamma $-symmetric structure 
on the homogeneous space $M=G/H$.
\end{proposition}


\subsection{Canonical connections on the homogeneous space $G/H$}\label{sscchs}


Let $(G/H,\Gamma)$ be a $\Gamma $-symmetric space. As we learned, the compexified Lie algebra $\g\otimes\C$ of $G$ is then $P$-graded, $P=\wg$, 
$\g\otimes\C=\oplus _{ p \in P}\  (\g\otimes\C)_ p$. If $\ve$ is the identity element of $P$ then the component $\h=(\g\otimes\C) _ 1$ 
is a Lie subalgebra of $\g\otimes\C$ of the form $\h\otimes\C$ where $\h=\g^\gm$, the set of fixed points of the action of $\gm$ on $\g$ and also the Lie algebra of the subgroup $H$. Let us consider the subspace $\g^\prime$ of $\g$:
$$\g^\prime=\oplus _{ p\neq  1} \ \g_ p .$$
For every $ 1\neq p\in P$, if $ p^2= 1$ then  $(\g\otimes\C)_ p=\g_ p\otimes\C$ where $\g_ p$ is given by (\ref{egc}), and if not, then $(\g\otimes\C)_ p\oplus(\g\otimes\C)_{ p^{-1}}=\widetilde{\g}_{ p}\otimes\C$ where $\widetilde{\g}_ p$ is the subspace of $\g$ spanned by the real and imaginary parts of the vectors in $(\g\otimes\C)_{ p}$. 

This simple claim follows because we have $\overline{\gamma(u+vi)}=\gamma(\overline{u+vi})$ where the action of $\gm$ on $\g\otimes\C$ is given by $\gamma(u+vi)=\gamma(u)+\gamma(v)i$. Thus if $ p^2= 1$ and $u+vi\in(\g\otimes\C)_{ p}$ then $\overline{\gamma(u+vi)}=\gamma( p)(u-vi)$ and $\gamma(\overline{u+vi})=\gamma(u)-\gamma(v)i$. Since $\gamma( p)$ is real, $\gamma(u)=\gamma( p)(u)$ and $\gamma(v)=\gamma( p)(v)$, proving that $u,v\in\g_{ p}$ and $(\g\otimes\C)_{ p}=\g_{ p}\otimes\C$. But if $ p^2\neq 1$  and again $u+vi\in(\g\otimes\C)_{ p}$ then $$\overline{\gamma(u+vi)}=\overline{\gamma( p)}(u-vi)=\gamma( p^{-1})\overline{u+vi}=\gamma(u-vi)$$
showing that the complex conjugation leaves invariant $(\g\otimes\C)_ p\oplus(\g\otimes\C)_{ p^{-1}}$. Then, of course, our claim follows.

If we set $\m$ the sum of all $\g_{ p}$ and $\widetilde{\g}_ p$ if $p\ne 1$  then $\g\otimes\C=\g_ 1\otimes\C\oplus\m\otimes\C$ and $\g=\g_ 1\oplus \m$. We also have
$$[\g_ 1,\m]\subset \m$$ 
so that $\m$ is an $\ad\g_ 1$-invariant subspace. If $H$ is a connected Lie group then $[\g_ 1,\m]\subset \m$ is equivalent to
$(\ad H)(\m)\subset \m$, that is, $\m$ is an $\ad H$-invariant subspace. This property is true without any conditions on $H$.

\begin{lemma}\label{lred}
Any $\gm$-symmetric space $(G/H,\Gamma)$ is reductive.
\end{lemma}

\begin{proof} Let us consider the associated local $\gm$-symmetric space $(\g/\h,\gm)$. We need to find a decomposition $\g=\h\oplus\m$ such that $\m$ is invariant under the adjoint action of the isotropy subgroup $H$ or, which is the same, under the action of the isotropy subalgebra $\h$. Now since $\h\otimes\C=(\g\otimes\C)_ 1$ and $\m\otimes\C=\oplus_{ p\neq 1}(\g\otimes\C)_{ p}$ we have that $[\h\otimes\C,\m\otimes\C]\subset\m\otimes\C]$. Then, of course, also $[\h,\m]\subset\m$, and the proof is complete.
\end{proof}

\medskip

We now deduce  from \cite{KN}, Chapter X, that $M=G/H$ admits  two $G$-invariant canonical connections 
denoted by $\nabla$ and $\overline \nabla$. The \emph{first canonical connection}, $\nabla$, satisfies
$$
\left\{
\begin{array}{l}
R(X,Y)=-\ad([X,Y]_\h), \ T(X,Y)_{\bar e}=-[X,Y]_{\m}, \ \forall X,Y \in \m \\
\nabla T=0 \\
\nabla R=0
\end{array}
\right.
$$
where $T$ and $R$ are the torsion and the curvature tensors of $\nabla$. The tensor $T$ is trivial
if and only if $[X,Y]_{\m}=0$ for all $X,Y \in \m$. This means that $[X,Y] \in \h$ that is $[\m,\m]\subset \h$.
Then the Lie algebra $\g$ is $\Z_2$-graded and the homogeneous space $G/H$ is symmetric. If the grading of $\g$ is given by $\Gamma$ where $\Gamma $ is not isomorphic to $\Z_2$, then $[\m,\m]$ need not be a subset of in $\h$ and then the torsion $T$ need not vanish.
In this case another connection, $\bar \nabla$, will be defined if one sets $\bar \nabla _X Y= \nabla_XY-T(X,Y)$. This is an affine invariant torsion free connection
on $G/H$ which has the same geodesics as $\nabla$. This connection is called the \emph{second canonical connection} or the \emph{torsion-free canonical
connection}. For example, if $\Gamma =\z$ then these connections can be distinct, as one can see from a number of examples in Section \ref{scl}.
\medskip

\noindent{\bf Remark. } Actually, there is another way of writing the canonical 
affine connection of a $\gm$-symmetric space, without any reference to Lie algebras. 
This is done by an intrinsic construction of $\Gamma$ - symmetric spaces  proposed by Lutz 
in \cite{L}.

\section{Classification of  local $\z$ - symmetric complex spaces}\label{scl}


We have seen that the classification of $\Gamma $-symmetric spaces $(G/H,\Gamma)$,
when $G$ is connected and simply connected,  corresponds  to the classification of Lie algebras graded by $P=\wg $. Below we establish the classification of  local $\z$-symmetric spaces $(\g,\gm)$ in the case where the corresponding Lie algebra $\g$ is simple complex and classical. 

We recall Definition \ref{deq} that given two $P$-gradings $\g=\oplus_{p\in P}\g_{p}$ and 
$\g=\oplus_{p\in P}\g_{p}^{\prime}$ of an algebra $\g$ by a group $P$ 
we call them \textit{equivalent} if there exists an automorphism $\alpha$ of $\g$ such that 
$\g_{p}^{\prime}=\alpha(\g_{p})$. To make the classification even more compact we will 
use another equivalence relation on the gradings. We will call
two $P$-gradings $\g=\oplus_{p\in P}\g_{p}$ and 
$\g=\oplus_{p\in P}\g_{a}^{\prime}$ of an algebra $\g$ by a group $P$  
\textit{weakly equivalent} if there exists an automorphism 
$\pi$ of $\g$ and an automorphism $\omega$ of $P$ such that 
$\g_{p}^{\prime}=\pi(\g_{\omega(p)})$. So the classification we are about to produce 
will be \textit{up to the weak equivalence}. 

\subsection{Introductory remarks about the $\z$-gradings}\label{ssir}

In this section $P =\{e,a,b,c\}$ is the group $\z$ with identity $e$ and $a^2=b^2=c^2=e, \ ab=c.$ We
will consider the $\z$-grading on a complex simple Lie algebra $\g$ one if the types $A_l, \ l\geq 1$, 
$B_l, \ l\geq 2$, $C_l, \ l \geq 3$ and $D_l, \ l\geq 4$.  We are going to use some results of    
 \cite{BShZ} and \cite{BZA}. That is why in the remainder of the paper we denote by $e$ the identity of the grading group $P$. 
 
Note that in those papers we do not consider the case of $\soo(8)$. All fine group gradings on this algebra have been described in \cite{DV}. We are grateful to the referee for pointing out that one can use this classification to construct all $\z$-gradings on $D_4$ by  Proposition 2 of  \cite{DM}. At the same time, we note that if one considers the $P$-gradings of $\g=\soo(8)$ where $P$ is an elementary abelian $2$-group, then each grading is equivalent to a grading induced from a grading of the matrix algebra $M_8$. This quickly follows from the description of the outer automorphisms of order 2. If we fix a canonical realization of $D_4$ in $M_8$ and a basis of the root system for $D_4$ then one of the three diagram automorphisms of order 2 can be given as the conjugation by an appropriate nonsingular matrix in $M_8$ while the others are the conjugates of this fixed one by the diagram automorphisms of order 3 (see, for example, \cite[Chapter III]{NJ}).

\medskip

\subsubsection{}\label{5.1.1} According to \cite{BShZ} any $P$-grading of a simple Lie algebra
$\g=\soo(2l+1), \ l\geq 2, \ \g=\soo(2l),\ l> 4$ and $\spp(2l),\ l\geq 3$ is induced from
an $P$-grading of the respective associative matrix algebra $R=M_{2l+1}$ in the first case, or
$M_{2l}$ in the second and the third case. As we just explained, this is also true for $\soo(8)$ provided that $P=\z$. Two kinds of $P$-grading on
the associative matrix algebra $M_n = R=\oplus_{ p\in P} \ R_a$ are of special importance:

1) \textit{Elementary gradings}. Each elementary grading is defined by an $n$-tuple  $(p_1,...,p_n)$ 
of elements of $P$ in such a way that all matrix units
$E_{ij}$ are homogeneous with $E_{ij} \in R_p $ if and only if $p =p_i^{-1}p_j.$

2) \textit{Fine gradings}. The characteristic property of such gradings is that for every 
$p  \in \su(R)$, $\dim R_p =1$ where $\su(R)=\{p\in P ,
\ {\mbox{\rm dim}}R_p \neq 0\}$. In the case of $P=\z$, each fine grading is either trivial or weakly equivalent 
to the grading on  $R\cong M_2$ 
given by the Pauli matrices $$X_e=I= \left( \begin{array}{ll} 1&0 \\ 0& 1 \end{array} \right),\ 
X_a= \left( \begin{array}{rr} -1&0 \\ 0& 1 \end{array} \right), \
X_b=\left( \begin{array}{ll} 0&1 \\ 1& 0 \end{array} \right), \
X_c=\left( \begin{array}{rr} 0&-1 \\ 1& 0 \end{array} \right)$$
in such a way that the graded component of degree $p$ is spanned by $X_{p}$, $p=e,a,b,c$. 
Notice \cite{BSZ} that the support of a fine grading of a simple associative algebra is always a subgroup of $P$.

According to \cite {BSZ} and \cite {surgrad} any $P$-grading of $R=M_n$ can be written as the tensor product of two 
graded matrix subalgebras $R=A\otimes B$, where $A\cong M_k$ and its grading is (equivalent to) elementary, and the
grading of $B\cong M_l$ is fine with $\su{A}\cap\: \su{B}=\{e\}$, $kl=n$. Thus the only cases possible, when $P=\z$, 
are 

1) $B=\C$ and $R=A\otimes \C=A$

2) $B=M_2$ and the grading on $A$ is trivial.

If $R$ is graded by $P$ as above, then an involution $\ast$ of $R=M_n$ is called \textit{graded} if 
$(R_{p})^{\ast}=R_{p}$ for any $p\in P$. In the case of such involution, the spaces 
$K(R,\ast )=\{X \in R, X^\ast =-X\}$ of skew - symmetric elements under $\ast$ and $H(R,\ast )=\{X \in R, X^\ast =X\}$ 
of symmetric elements under $\ast$ are graded and the first is a simple Lie algebra of one of the types $B, C, D$ under
the bracket $[X,Y]=XY-YX$.

It is proved in \cite{BShZ} that $\g$ as a $P$ - graded algebra is isomorphic to $K(R,\ast )=\{X \in R, X^\ast = -X\}$ 
for an appropriate graded involution $\ast$. 

In general the involution does not need to respect $A$ and $B$. 
But this is definitely the case when $P =\z$. In fact, since the grading of $B$ is fine, 
using the support conditions, we see that either $B=\C$ and $R=A\otimes \C=A$ is respected by $\ast $, 
or $B=M_2$ and the grading of $A$ is trivial.  Since $B$ is the centralizer of $A$ in $R$ it follows that $B$
is also invariant under the involution. For the details of the above claims see \cite{BZA}.

Now any involution has the form $\ast: X\longrightarrow X^\ast =\Phi ^{-1}\tr{X}\Phi$, for a nonsingular matrix $\Phi$, 
which is either symmetric in 
the orthogonal case and skew-symmetric in the symplectic case. Since the elementary and fine 
components are invariant under the involution, we have that $\Phi=\Phi_1\otimes \Phi_2$ where 
$\Phi_1$ defines a graded involution on $A$ and $\Phi_2$ on $B$.

First we recall the description of the graded involutions for the elementary gradings. Given an element $p$ of a group $P$ and a natural number $k$ we denote by $p^{(n)}$ the $n$-tuple $\underbrace{p,\ldots,p}_{k}$.
Using the argument of \cite{BZA}, one may assume that this grading is
given by an $n$-tuple $\nu=( p _1^{(k_1)},p _2^{(k_2)},...p _{s+2t}^{(k_{s+2t})})$, $n=k_1+..+k_{s+2t}$, where
$p _i \in \Gamma $ are pairwise different, $k_{s+1}=k_{s+2}, \ldots, k_{s+2t-1}=k_{s+2t}$, 
and there is $p_0$ 
such that $p _0=p_1^2=...=p_s^2=
p _{s+1}p _{s+2}=...=p _{s+2t-1}p _{s+2t}.$ However, since $p ^2=e$ for the elements of $\z$,
we have $p _0=e$ and if $st \neq 0$ we have $p _{s+1}p _{s+2}=e$ this implies $p _{s+1}=p _{s+2}$
which gives a contradiction. Then in this case either $t=0$ and the $\z$-elementary grading corresponds to
$e=a^2=b^2= c^2$ and is given by the $p$-tuple
$$(e^{(k_1)},a^{(k_2)},b^{(k_3)},c^{(k_4)})$$ 
or $s=0$ and the $\z$-grading  corresponds to $a= ae=bc$ with $k_{1}=k_2$ and $k_3=k_4$. 

 In the general case the matrix $\Phi $ defining the involution has the form
$$
\Phi=\mathrm{diag}\left\{I_{k_1},\ldots,I_{k_s},\left(\begin{array}{cc}0&I_{k_{s+1}}\\I_{k_{s+1}}&0\end{array}\right),
\ldots,\left(\begin{array}{cc}0&I_{k_{s+2t-1}}\\I_{k_{s+2t-1}}&0\end{array}\right)\right\}
$$
if $\ast $ is a transpose involution, i.e. $\Phi$ is symmetric, or
$$
\Phi=\mathrm{diag}\left\{S_{k_1},\ldots,S_{k_p},\left(\begin{array}{cc}0&I_{k_{p+1}}\\-I_{k_{p+1}}&0\end{array}\right),\ldots,\left(\begin{array}{cc}0&I_{k_{p+2s-1}}\\-I_{k_{p+2s-1}}&0\end{array}\right)\right\}.
$$
in the case of a skew-symmetric $\Phi$, where we denote by $I_k$ the identity matrix of order $k$ and by $S_{2l}$ the skew-symmetric 
matrix 
$\left(
\begin{array}{cc}
0&I_l\\
-I_l&0
\end{array}
\right)$. 
If $\Gamma =\z$, when we consider the case $t=0$, the matrices of $\Phi$ are the identity in the symmetric case and
$$
\Phi=\mathrm{diag}\left\{S_{k_1},\ldots,S_{k_4}\right\}
$$ in the skew-symmetric case and if $s=0$ then
$$
\Phi=\mathrm{diag}\left\{
\left(
\begin{array}{cc}
0&I_{k_{1}}
\\I_{k_{1}}&0
\end{array}
\right),
\left(\begin{array}{cc}0&I_{k_{3}}\\I_{k_{3}}&0\end{array}\right)\right\}
$$
if $\Phi$ is symmetric and
$$
\Phi=\mathrm{diag}\left\{\left(\begin{array}{cc}0&I_{k_{1}}\\-I_{k_{1}}&0\end{array}\right),
\left(\begin{array}{cc}0&I_{k_{3}}\\-I_{k_{3}}&0\end{array}\right)\right\}
$$
if $\Phi$ is skew-symmetric.

\medskip

 If $R=A\otimes B$ and $B \neq \C$, then $\su(A) = \{e\}$. We have $\Phi=\Phi _1\otimes \Phi _2$
and the involution on $R$ defines involutions on $A$ and $B$. It follows that $\Phi $ is symmetric if and only if
either $\Phi _1$ and $\Phi _2$ are both symmetric or they both are skew-symmetric. Similarly, $\Phi $ is 
skew-symmetric if one of $\Phi _1$,$\Phi _2$ is symmetric and the other is skew-symmetric. It was proved
in  \cite{BShZ} that $M_2$ with graded involution is isomorphic to $M_2$ with $\Gamma $-graded basis
$$
\left\{X_e=I_2,X_a=\left( \begin{array}{cc}
-1 & 0 \\
0 & 1
\end{array}
\right),
X_b=\left( \begin{array}{cc}
0 & 1 \\
1 & 0
\end{array}
\right),
X_c=\left( \begin{array}{cc}
0 & -1 \\
1 & 0
\end{array}
\right)
\right\}
$$
and the graded involution is given by one of $\Phi =X_e,X_a,X_b,X_c.$

\bigskip

\subsubsection{}\label{5.1.2} Now let $\g$ be a simple Lie algebra of type $A_l$. We view $\g$ as the set $\g=\sll(n)$  of all matrices of trace zero in 
the matrix algebra $R=M_n(\C)$ where $n=l+1$. In this case any $P$-grading of $\g$ belongs to one of the following two classes (see \cite{BZA}).

- For Class $\mathrm{I}$ gradings, any grading of $\g$ is induced from a $P$-grading of $R=\bigoplus_{p\in P}R_{p}$ and one simply has to set 
$\g_p =R_p $ for $p  \neq e$ and $\g_e=R_e\cap \g$, otherwise. For $P =\z$ we still have to
distinguish between the cases $R=A\otimes\C$ with an elementary grading on $A=M_n$ or $R=A\otimes B$ with trivial
grading on $A$ and fine $\z$-grading on $B=M_2$.

- For Class $\mathrm{II}$ gradings, we have to fix an element $q$ of order $2$ in $P$ and an involution 
$P$-grading $R=\oplus _{p  \in P } \  R_p $. Then for any $p  \in P $ one has
$$
\g_p =K(R_p ,\ast )\oplus H(R_{p q},\ast )\cap \g.
$$
The involution grading on $M_n$ have been discussed just before in 4.1.1. It should be noted that in the case where 
$B \neq \C$
in $ R=A \otimes  B$ we have
$$(\ast ) \ \ K(R_p ,\ast )=K(A,\ast )\otimes H(B_p ,\ast )
\oplus H(A,\ast )\otimes K(B_p ,\ast ),$$
$$(\ast \ast  ) \ \ H(R_p ,\ast )=H(A,\ast )\otimes H(B_p ,\ast )\oplus 
K(A,\ast )\otimes K(B_p ,\ast ).$$
As noted above, $\Phi =\Phi _1\otimes \Phi _2$ with $\Phi _2=X_p $ for $p =e,a,b,c$. If $X_p $
is symmetric with respect to $\Phi _2$ then $(\ast )$ and $(\ast \ast )$ become
$$ K(R_p ,\ast )=K(A,\ast )\otimes X_p ,$$
$$  H(R_p ,\ast )=H(A,\ast )\otimes X_p .$$
If $X_p $ is skew-symmetric with respect to $\Phi _2$ then
$$ K(R_p ,\ast )=H(A,\ast )\otimes X_p ,$$
$$  H(R_p ,\ast )=K(A,\ast )\otimes X_p .$$

All these remarks allow to determine up to the weak equivalence the pairs $(\g,\g_e)$ inside the respective matrix algebra $M_n$.
This gives a local classification of $\z$-symmetric homogeneous spaces $G/H$ where $G$ is simple classical connected Lie group.

\subsection{Classification : $B_l,C_l,D_l$ cases}\label{sscbcd}

We have that $\g=K(M_n,\Phi )$, $\Phi $ being symmetric for the cases $B_l$ and $D_l$ and $\Phi $ skew-symmetric
for $C_l$. In the case of $B_l$ we have $n=2l+1$ and in the case of $C_l$ and $D_l$ we have $n=2l$. 

\subsubsection{Lie gradings corresponding to the elementary grading of $M_n$}\label{5.2.1} 

Since we are interested in the gradings only up to the weak equivalence, it is sufficient to consider the following tuples defining the elementary gradings:
$$
\left\{
\begin{array}{l}
\nu _1=(e^{(k_1)},a^{(k_2)})\\
\nu _2=(e^{(k_1)},a^{(k_2)},b^{(k_3)})\\
\nu _3=(e^{(k_1)},a^{(k_2)},b^{(k_3)},c^{(k_4)}).
\end{array} 
\right.
$$
Notice that the $\z$-gradings corresponding to $\nu_1$ coincide with $\Z_2$-gradings and thus the corresponding homogeneous spaces are  
symmetric in the classical sense. The matrices defining the graded transpose involution in the case of $\nu_1$ are
$$
\Phi _1=\mathrm{diag}\left\{I_{k_1},I_{k_2}\right\}\mathrm{ and  }\; \Phi _1'=\left(\begin{array}{cc}0&I_{k_{1}}\\I_{k_{1}}&0\end{array}\right).
$$
In the case of $\nu_2$ we have
$$\Phi _2=\mathrm{diag}\left\{I_{k_1},I_{k_2},I_{k_3}\right\}.$$
Finally, in the case of  for $\nu_3$ we have
$$\Phi _3=\mathrm{diag}\left\{I_{k_1},I_{k_2},I_{k_3},I_{k_4}\right\} \ {\mbox{\rm or}} \
\Phi _3'=\mathrm{diag}\left\{\left(\begin{array}{cc}0&I_{k_{1}}\\I_{k_{1}}&0\end{array}\right),
\left(\begin{array}{cc}0&I_{k_{3}}\\I_{k_{3}}&0\end{array}\right)\right\}.
$$
 If the involution is symplectic, then the respective matrices, in the case of $\nu_1$, are as follows:
$$
\overline\Phi _1=\mathrm{diag}\left\{S_{k_1},S_{k_2}\right\}\mathrm{ and }\; \overline\Phi _1'=\left(\begin{array}{cc}0&I_{k_{1}}\\-I_{k_{1}}&0\end{array}\right).
$$
In the case of $\nu_2$:
$$\overline\Phi _2=\mathrm{diag}\left\{S_{k_1},S_{k_2},S_{k_3}\right\}$$
 and in the case of $\nu_3$:
$$\overline\Phi _3=\mathrm{diag}\left\{S_{k_1},S_{k_2},S_{k_3},S_{k_4}\right\} \ {\mbox{\rm or}} \
\overline\Phi' _3=\mathrm{diag}\left\{\left(\begin{array}{cc}0&I_{k_{1}}\\-I_{k_{1}}&0\end{array}\right),
\left(\begin{array}{cc}0&I_{k_{3}}\\-I_{k_{3}}&0\end{array}\right)\right\}.
$$

In the cases $(\nu _1,\Phi _1)$ and $(\nu _2,\overline \Phi _1)$ we have $\g_e=\soo(k_1) \oplus \soo(k_2)$
and $\g_e=\spp(k_1) \oplus \spp(k_2)$, respectively. 

In the cases $(\nu _1,\Phi ' _1)$ and $(\nu _3,\overline \Phi ' _1)$ we have that
$$\g_e=\left\{\left.\left(\begin{array}{cc}U_1&0\\0&-\tr{U}_1\end{array}\right)\right|
 U_1 \in M_{k_1}\right\}.$$ So for both $D_{k_1}$ and $C_{k_1}$  
cases we will have $\g_e=\gll(k_1)$.

In the cases $(\nu _2,\Phi _2)$ and $(\nu _2,\overline \Phi _2)$ we have $\g_e=\soo(k_1) \oplus \soo(k_2) \oplus \soo(k_3)$
and $\g_e=\spp(k_1) \oplus \spp(k_2) \oplus \spp(k_3)$, respectively. 

In the cases $(\nu _3,\Phi _3)$ and $(\nu _3,\overline \Phi _3)$ we have 
$\g_e=\soo(k_1) \oplus \soo(k_2) \oplus \soo(k_3) \oplus \soo(k_4)$
and $\g_e=\spp(k_1) \oplus \spp(k_2) \oplus \spp(k_3)\oplus \spp(k_4)$, respectively.
 
In the cases $(\nu _3,\Phi ' _3)$ and $(\nu _3,\overline \Phi ' _3)$ we have that
$$\g_e=\mathrm{diag}\left\{\left(\begin{array}{cc}U_1&0\\0&-^tU_1\end{array}\right),
\left(\begin{array}{cc}U_2&0\\0&-^tU_2\end{array}\right)\right\}, \;
U_1 \in M_{k_1}\mathrm{ and }\; U_2 \in M_{k_2}.$$ So for both $D_{k_1+k_3}$ and $C_{k_1+k_3}$  
cases we have $\g_e=\gll(k_1) \oplus \gll(k_3)$.

\bigskip

\begin{center}
\begin{tabular}{|c|c|}
\hline
&\\
$\g$ & $\g_e$ \\
& \\
\hline\hline &\\
$\soo(k_1 +k_2)$ & $\soo(k_1) \oplus \soo(k_2)$ \\
& \\
\hline &\\
$\soo(k_1 +k_2 +k_3)$ & $\soo(k_1) \oplus \soo(k_2) \oplus \soo(k_3)$ \\
& \\
\hline &\\
$\soo(k_1 +k_2 +k_3+k_4)$ & $\soo(k_1) \oplus \soo(k_2) \oplus \soo(k_3) \oplus \soo(k_4)$ \\
& \\
\hline &\\
$\spp(k_1 +k_2)$ & $\spp(k_1) \oplus \spp(k_2)$ \\
& \\
\hline &\\
$\spp(k_1 +k_2 +k_3)$ & $\spp(k_1) \oplus \spp(k_2) \oplus \spp(k_3)$ \\
& \\
\hline &\\
$\spp(k_1 +k_2 +k_3+k_4)$ & $\spp(k_1) \oplus \spp(k_2) \oplus \spp(k_3) \oplus \spp(k_4)$ \\
& \\
\hline &\\
$\soo(2m)$ & $\gll(m)$ \\
& \\
\hline &\\
$\soo(2(k_1+k_2))$ & $\gll(k_1) \oplus \gll(k_2)$ \\
& \\
\hline &\\
$\spp(2m)$ & $\gll(m)$ \\
& \\
\hline &\\
$\spp(2(k_1+k_2))$ & $\gll(k_1) \oplus \gll(k_2)$\\& \\
\hline
\end{tabular}

\bigskip

\textit{Table 1}
\end{center}
{\bf Conclusion.} The $\z$-symmetric spaces corresponding to the pairs
$(\g,\g_e)$ considered so far are given in the  Table 1.

In all the above cases the $\z$-gradings of $\g$
are very easy to compute, using explicit matrices of the involutions. For example, the gradings of $\soo(n)$ are given in \cite{BGR}.

\medskip

\subsubsection{Gradings corresponding to  $R=A \otimes B$, with $B$ nontrivial}\label{sssbnon}

The algebra $B$ is endowed with a fine grading given by the Pauli matrices
$$X_e=I_2,\ X_a= \left( \begin{array}{ll} -1&0 \\ 0& 1 \end{array} \right), \
X_b=\left( \begin{array}{ll} 0&1 \\ 1& 0 \end{array} \right), \
X_c=\left( \begin{array}{ll} 0&-1 \\ 1& 0 \end{array} \right)$$
and  $A=M_m$, $\Phi =\Phi_1 \otimes  \Phi_2$. We have that
$$\g=K(R,\Phi )=K(A,\Phi _1) \otimes H(B,\Phi _2) \oplus H(A,\Phi _1) \otimes K(B,\Phi _2).$$
In particular, $\g_e =K(R_e,\Phi _1)\otimes I_2$. If $\Phi _1$ is symmetric then $\g_e =\soo(m)$, and
if $\Phi _1$ is skew symmetric, $g_e=\spp(m).$

\medskip

 {\bf Conclusion.} We obtain the $\z$-symmetric spaces corresponding to the pairs
$(\g,\g_e)$ given in the following table:

\begin{center}
\begin{tabular}{|c|c|}
\hline
&\\
$\g$ & $\g_e$ \\
& \\
\hline\hline
& \\
$\soo(2m)$ & $\soo(m)$ \\
& \\
\hline
& \\
$\soo(4m)$ & $\spp(2m)$ \\
& \\
\hline
& \\
$\spp(4m)$ & $\spp(2m)$\\
& \\
\hline
& \\
$\spp(2m)$ & $\soo(m)$ \\  
& \\
\hline
\end{tabular}

\bigskip

\textit{Table 2}
\end{center}

\medskip

In these cases we will describe the components of the gradings explicitly. 

If $\g=\soo(2m)$ then $\Phi $ is symmetric
and is of one of the form
$$ \Psi _1=I_m \otimes I_2, \ \Psi _2=I_m \otimes X_a, \ \Psi _3=I_m \otimes X_b, \ \Psi _4=S_m \otimes X_c.   $$
For $\Psi  _1$ we have 
$$\g=\g(\Psi_1)=K(A,I_m)\otimes I\oplus K(A,I_m)\otimes X_a\oplus K(A,I_m)\otimes X_b \oplus H(A,I_m) \otimes X_c,$$ 
for $\Psi  _2$ we have $$\g=\g(\Psi_2)=K(A,I_m)\otimes I\oplus K(A,I_m)\otimes X_a\oplus H(A,I_m) \otimes X_b \oplus K(A,I_m)\otimes X_c,$$
for $\Psi  _3$ we have $$\g=\g(\Psi_3)=K(A,I_m)\otimes I \oplus H(A,I_m) \otimes X_a\oplus K(A,I_m)\otimes X_b\oplus K(A,I_m)\otimes X_c,$$ and
for $\Psi  _4$ we have 
$$\g=\g(\Psi_4)=K(A,S_m)\otimes I \oplus H(A,S_m) \otimes X_a\oplus H(A,S_m) \otimes X_b\oplus H(A,S_m) \otimes X_c.$$
The conjugation by $I_m\otimes\left( \begin{array}{ll} i&0 \\ 0& 1 \end{array} \right)$ maps 
$\g(\Psi_1)$ to $\g(\Psi_2)$ while mapping $K(A,I_m)\otimes I$ and $K(A,I_m)\otimes X_a$ into themselves, $K(A,I_m)\otimes X_b$ into $K(A,I_m)\otimes X_c$ and $H(A,I_m)\otimes X_c$ into $H(A,I_m)\otimes X_b$. If we apply an automorphism of $\z$ changing places of $b$ and $c$ we will see that the first and the second gradings are weakly equivalent. Quite similarly, the conjugation by $I_m\otimes\frac{1}{\sqrt{2}}\left( \begin{array}{lr} 1&1 \\ i& -i \end{array} \right)$ maps 
$\g(\Psi_1)$ to $\g(\Psi_3)$ while mapping $K(A,I_m)\otimes I$ into itself, $K(A,I_m)\otimes X_a$ into $K(A,I_m)\otimes X_b$, $K(A,I_m)\otimes X_b$ into $K(A,I_m)\otimes X_c$ and $H(A,I_m)\otimes X_c$ into $H(A,I_m)\otimes X_a$. It remains to apply an automorphism of $\z$ mapping $a$ to $c$, $c$ to $b$ and $b$ to $a$ to make sure that the first and the third gradings are weakly equivalent. Thus, all the first three gradings are weakly equivalent. None of them is weakly equivalent to the fourth one because in these cases  $\g_e\cong\soo(m)$ while in the fourth case we have  $\g_e\cong\spp(m)$.

The matrix form of the first and the fourth gradings are given below. 

For $\Psi  _1$ we have
$$\g= \left\{ \left( \begin{array}{ll} U_1-U_2&U_3-V \\ U_3+V& U_1+U_2 \end{array} \right), \ U_1,U_2,U_3 \in \soo(m), \
^tV=V \right\} .$$ The components are:
$$\g_e \oplus \g_a \oplus \g_b \oplus g_c = 
\left\{ \left( \begin{array}{ll} U_1&0 \\ 0& U_1 \end{array} \right) \right\} \oplus 
\left\{ \left( \begin{array}{ll} -U_2&0 \\ 0& U_2 \end{array} \right)\right\} \oplus
\left\{ \left( \begin{array}{ll} 0&U_3\\ U_3& 0 \end{array} \right) \right\} \oplus
\left\{ \left( \begin{array}{ll} 0&-V \\ V& 0 \end{array} \right) \right\}. $$   

For $\Psi  _4$ we have
$$\g= \left\{ \left( \begin{array}{ll} P-Q_1&Q_2-Q_3 \\ Q_2-Q_3& P+Q_1 \end{array} \right), \ P \in \spp(m), \ 
Q_1,Q_2,Q_3 \in H(M_m,S_m) \right\} .$$ The components are:
$$\g_e \oplus \g_a \oplus \g_b \oplus g_c = 
\left\{ \left( \begin{array}{ll} P&0 \\ 0& P \end{array} \right) \right\} \oplus 
\left\{ \left( \begin{array}{ll} -Q_1&0 \\ 0& Q_1 \end{array} \right)\right\} \oplus
\left\{ \left( \begin{array}{ll} 0&Q_2\\ Q_2& 0 \end{array} \right) \right\} \oplus
\left\{ \left( \begin{array}{ll} 0&-Q_3 \\ -Q_3& 0 \end{array} \right) \right\}. $$
 
If $\g=\spp(2m)$ then $\Phi $ is skew-symmetric
and is of one of the form
$$ \overline \Psi _1=S_m \otimes I_2, \ \overline \Psi _2=S_m \otimes X_a, 
\ \overline \Psi _3=S_m \otimes X_b, \ \overline \Psi _4=I_m \otimes X_c.   $$
For $\overline \Psi _1$ we have 
$$\g=\g(\overline \Psi _1)=K(A,S_m)\otimes I\oplus K(A,S_m)\otimes X_a\oplus K(A,S_m)\otimes X_b \oplus H(A,S_m) \otimes X_c,$$ 
for $\overline \Psi _2$ we have $$\g=\g(\overline \Psi _2)=K(A,S_m)\otimes I\oplus K(A,S_m)\otimes X_a\oplus H(A,S_m) \otimes X_b \oplus K(A,S_m)\otimes X_c,$$
for $\overline \Psi _3$ we have $$\g=\g(\overline \Psi _3)=K(A,S_m)\otimes I \oplus H(A,I_m) \otimes X_a\oplus K(A,I_m)\otimes X_b\oplus K(A,I_m)\otimes X_c,$$ and
for $\Psi  _4$ we have 
$$\g=\g(\Psi_4)=K(A,S_m)\otimes I \oplus H(A,S_m) \otimes X_a\oplus H(A,S_m) \otimes X_b\oplus H(A,S_m) \otimes X_c.$$

The same argument as before shows that the first three gradings are weakly equivalent and none of them is weakly equivalent to the fourth one. In the first three cases we have $\g_e\cong\spp(m)$ while in the fourth case we have  $\g_e\cong\soo(m)$.

Again we give the matrix form of the first and the fourth gradings. 
If $\Phi $ is skew symmetric then $\g=\spp(2m)$ and $\Phi $ is of one of the form
$$ \overline \Psi _1=S_m \otimes I_2, \ \overline \Psi _2=S_m \otimes X_a, 
\ \overline \Psi _3=S_m \otimes X_b, \ \overline \Psi _4=I_m \otimes X_c.   $$
For $\overline \Psi  _1$ we have 
$$\g= \left\{ \left( \begin{array}{ll} U_1-U_2&U_3-V \\ U_3+V& U_1+U_2 \end{array} \right), \ U_1,U_2,U_3 \in \spp(m), \
V \in H(A,S_m) \right\} .$$ The components are:
$$\g_e \oplus \g_a \oplus \g_b \oplus g_c = 
\left\{ \left( \begin{array}{ll} U_1&0 \\ 0& U_1 \end{array} \right) \right\} \oplus 
\left\{ \left( \begin{array}{ll} -U_2&0 \\ 0& U_2 \end{array} \right)\right\} \oplus
\left\{ \left( \begin{array}{ll} 0&U_3\\ U_3& 0 \end{array} \right) \right\} \oplus
\left\{ \left( \begin{array}{ll} 0&-V \\ V& 0 \end{array} \right) \right\}. $$   
For $\overline \Psi  _4$ we have 
$$\g= \left\{ \left( \begin{array}{ll} P-Q_1&Q_2-Q_3 \\ Q_2-Q_3& P+Q_1 \end{array} \right), \ P \in \soo(m), \ 
Q_1,Q_2,Q_3 \in H(M_m,I_m) \right\} .$$ The components are:
$$\g_e \oplus \g_a \oplus \g_b \oplus g_c = 
\left\{ \left( \begin{array}{ll} P&0 \\ 0& P \end{array} \right) \right\} \oplus 
\left\{ \left( \begin{array}{ll} -Q_1&0 \\ 0& Q_1 \end{array} \right)\right\} \oplus
\left\{ \left( \begin{array}{ll} 0&Q_2\\ Q_2& 0 \end{array} \right) \right\} \oplus
\left\{ \left( \begin{array}{ll} 0&-Q_3 \\ Q_3& 0 \end{array} \right) \right\}. $$   
  
\subsection{Classification of Class $\mathrm{I}$ gradings on $A_l$-type Lie algebras}\label{ssA}
  
If no fine component is present in $R=M_n\supset \g=\sll(n)$, $n=l+1$, then all is defined by the $n$-tuples 
$$
\left\{
\begin{array}{l}
\nu _1=(e^{(k_1)},a^{(k_2)})\\
\nu _2=(e^{(k_1)},a^{(k_2)},b^{(k_3)})\\
\nu _3=(e^{(k_1)},a^{(k_2)},b^{(k_3)},c^{(k_4)}).
\end{array} 
\right.
$$
In the case of $\nu _1$, the grading correspond to a symmetric decomposition (in fact we obtain the symmetric pair
$$(\sll(n),\sll(k_1)\oplus \sll(k_2)\oplus\C)$$
(or $\R$ if we are in the real case.)

 In the case of $\nu_2$ the
$$\g_e=\mathrm{diag}\left\{\left. X,Y,Z, \right| X \in M_{k_1},Y \in M_{k_2},Z \in M_{k_3}, tr(X+Y+Z)=0\right\}$$
and $\g_e =\sll(k_1)\oplus \sll(k_2)\oplus \sll(k_3)\oplus \C^2.$

 In the case of $\nu_3$ the
$$\g_e=\mathrm{diag}\left\{\left. X,Y,Z,T \right| X \in M_{k_1},Y \in M_{k_2},Z \in M_{k_3},T \in M_{k_4},tr(X+Y+Z+T)=0 \right\}$$
and $\g_e =\sll(k_1)\oplus \sll(k_2)\oplus \sll(k_3)\oplus \sll(k_4)\oplus \C^3.$

 In all these case the grading is obvious. If $R=A\otimes B=M_{n}$, $A=M_m$, $n=2m$, with a trivial grading on $A$, then 
$$\g_e= \left\{\left.\left(\begin{array}{ll}X & 0 \\ 0 & X\end{array}\right)\right| X \in sl(m)\right\}$$
and the grading is given by
\begin{eqnarray*}\g&=&\g_e \oplus \g_a \oplus \g_b \oplus g_c\\ &=& \g_e \oplus \left\{\left.\left(\begin{array}{rr}X & 0 \\ 0 & -X\end{array}\right) \right| X \in M_m\right\}
\oplus \left\{\left.\left(\begin{array}{rr}0 & X \\ X & 0\end{array}\right)\right| X \in M_m\right\}\\
&\oplus& \left\{\left.\left(\begin{array}{rr}0 & -X \\ X & 0\end{array}\right) \right| X \in M_m)\right\}.\end{eqnarray*}

 {\bf Conclusion.} We obtain the $\z$-symmetric spaces corresponding to the pairs
$(\g,\g_e)$ given in the following table:

\begin{center}
\begin{tabular}{|c|c|}
\hline
& \\
$\g$ & $\g_e$ \\
& \\
\hline\hline 
& \\
$\sll(2n)$ & $\sll(n)$\\
& \\
\hline
& \\
$\sll(k_1+k_2)$ & $\sll(k_1)\oplus \sll(k_2)\oplus  \C$ \\
&\\ \hline 
& \\
$\sll(k_1+k_2+k_3)$ & $\sll(k_1)\oplus \sll(k_2)\oplus \sll(k_3)\oplus  \C^2$\\
& \\
\hline 
& \\
$\sll(k_1+k_2+k_3+k_4)$ & $\sll(k_1)\oplus \sll(k_2)\oplus \sll(k_3)\oplus \sll(k_4)\oplus \C^3$\\
& \\ \hline
\end{tabular}

\bigskip

\textit{Table 3}
\end{center}

\bigskip

\subsection{Classification of Class $\mathrm{II}$ gradings on $A_l$-type Lie algebras}

The general approach described in \ref{5.1.2} enables one to classify the Class $\mathrm{II}$ gradings on $\g=\sll(n)$, for any $n\ge 2$ and any grading group $P$. However, in the case $P=\z$ the amount of work can be significantly reduced if one uses the results of \cite{BZ} and \cite{BZA}. In the former paper, in the case of outer gradings of $\sll(n)$, the authors showed that the dual $\gm$ of the grading group $P$ decomposes as the direct product $\langle\vp\rangle\times\Lambda$ where $\vp$ is an antiautomorphism of order $2^m$ and $\Lambda$ acts by inner matrix automorphisms. If $H=\Lambda^{\perp}$ then $H$ is a subgroup of order 2 and the induced $G/H$-grading of $\sll(n)$ is induced from $M_n$. To obtain the $G$-grading of $\g$ one has to refine the $G/H$-grading by intersecting them with the eigenspaces of $\vp$. In the case $P=\z$, $\vp$ is (a negative to) a graded involution of $M_n$, $\Lambda$ is generated by an automorphism $\lambda$ of order $2$, the generators of $P$ are $a$ and $b$ such that $\lambda(a)=-1$, $\lambda(b)=1$, $\vp(a)=1$, and $\vp(b)=-1$. In the latter paper the authors described all graded involutions on graded matrix algebras. In our particular case the gradings by $\z$ on $\sll(n)$ correspond to $\mathbb{Z}_2$-graded eigenspaces of a (negative to) a graded involution $\vp$ on a $\mathbb{Z}_2$-graded associative algebra $R=M_n$. Any $\Z_2$-grading on $R$ is elementary, given by a tuple $\nu _1=(e^{(k_1)},a^{(k_2)})$ where $a$ is the generator of $\Z_2$ and $k_1+k_2=n$. Now, as described in \cite[Theorem 3]{BZA}, any graded involution is graded equivalent to $\omega(X)=\cj{\Phi}{\tr{X}}$ where $\Phi$ is of one the following types
\begin{eqnarray}\label{egrinv}	
\Phi _1&=&\mathrm{diag}\left\{I_{k_1},I_{k_2}\right\}\\
\Phi _1'&=&\left(\begin{array}{cc}0&I_{k_{1}}\\I_{k_{1}}&0\end{array}\right)\\
\overline\Phi _1&=&\mathrm{diag}\left\{S_{k_1},S_{k_2}\right\}\\
\overline\Phi _1'&=&\left(\begin{array}{cc}0&I_{k_{1}}\\-I_{k_{1}}&0\end{array}\right).
\end{eqnarray}
Now it remains to apply \cite[Corollary 5.6]{BZ}, where $K=\langle a\rangle$ to obtain that all Class $\mathrm{II}$ gradings of $\g$ have the following form
\begin{eqnarray*}
\g_e&=&K(R_e,\Phi)\\
\g_a&=&K(R_a,\Phi)\\
\g_b&=&H(R_{e},\Phi)\\
\g_c&=&H(R_{a},\Phi)
\end{eqnarray*}

In all four cases with, depending on the choice of $\Phi$, we have $$R=\left\{\left.\left(\begin{array}{cc}U&V\\W&T
\end{array}\right)\right|U\in M_{k_1}, T\in M_{k_2}\right\},\mbox{ and }\g=\left\{\left.\left(\begin{array}{cc}U&V\\W&T
\end{array}\right)\right| \mathrm{tr}\, U + \mathrm{tr}\, T=0\right\}$$ and in the last two cases, we additionally have that $k_1=k_2$. It easily follows that for the $\Z$- grading we have
$R_e=\left\{\left(\begin{array}{cc}U&0\\0&T
\end{array}\right)\right\}\subset R$ and $R_a=\left\{\left(\begin{array}{cc}0&V\\W&0
\end{array}\right)\right\}\subset R$. 

Now we can explicitly write all the gradings in this case. In what follows we keep the following notation. For any $X$ of size $c\times d$, we denote by $X^{\ast}$ the ordinary transpose of $X$, except in the case of $\Phi_1^{\prime}$ where it means $\inv{S_p}\tr{X}S_q$. By $U, U_1,\ldots$ we will denote the matrices with $X^{\ast}=-X$, $U, U_1,\ldots$ those with $X^{\ast}=X$, while by $W,W_1,\ldots$ we will mean any matrices of appropriate sizes. All matrices must be in $\g$.
 
In both cases $\Phi=\Phi_1, \overline\Phi _1$ then we will have
$$\g_e \oplus \g_a \oplus \g_b \oplus g_c = 
\left\{ \left( \begin{array}{ll} U_1&0 \\ 0& U_2 \end{array} \right) \right\} \oplus 
\left\{ \left( \begin{array}{ll} 0&W \\ -\as{W}& 0 \end{array} \right)\right\} \oplus
\left\{ \left( \begin{array}{ll} V_1&0\\ 0&V_2 \end{array} \right) \right\} \oplus
\left\{ \left( \begin{array}{ll} 0&W \\ \as{W}& 0 \end{array} \right) \right\}. $$ 
However, because we have the transpose involution in the first case and the symplectic in the second, we obtain two inequivalent local symmetric spaces $$\sll(k_1+k_2)/(\mathrm{so}\,(k_1)\oplus\mathrm{so}\,(k_2))$$
and
$$\sll(k_1+k_2)/(\mathrm{sp}\,(k_1)\oplus\mathrm{sp}\,(k_2)).$$

It should be noted that $k_1$ is always nonzero while we could have $k_2=0$. In this case we actually have a $\Z_2$-grading rather than a $\z$-grading. The respective local symmetric spaces are well-known to be
$$\sll(n)/\mathrm{so}\,(n)$$
and
$$\sll(2m)/\mathrm{sp}\,(2m).$$

In the case of $\Phi=\Phi_1^{\prime}$ we will have
$$\g_e \oplus \g_a \oplus \g_b \oplus g_c = 
\left\{ \left( \begin{array}{ll} W&0 \\ 0& -\tr{W} \end{array} \right) \right\} \oplus 
\left\{ \left( \begin{array}{ll} 0&U_1 \\ U_2& 0 \end{array} \right)\right\} \oplus
\left\{ \left( \begin{array}{ll} W&0 \\ 0& \tr{W} \end{array} \right) \right\} \oplus
\left\{ \left( \begin{array}{ll} 0&V_1\\ V_2&0 \end{array} \right) \right\}. $$ 

Finally, in the case of $\Phi=\overline{\Phi}_1^{\prime}$ we will have
$$\g_e \oplus \g_a \oplus \g_b \oplus g_c = 
\left\{ \left( \begin{array}{ll} W&0 \\ 0& -\tr{W} \end{array} \right) \right\} \oplus 
\left\{ \left( \begin{array}{ll} 0&V_1 \\ V_2& 0 \end{array} \right)\right\} \oplus
\left\{ \left( \begin{array}{ll} W&0 \\ 0& \tr{W} \end{array} \right) \right\} \oplus
\left\{ \left( \begin{array}{ll} 0&U_1\\ U_2&0 \end{array} \right) \right\}. $$

Obviously, these two latter grading are weakly equivalent, and the weak equivalence is achieved by an automorphism of $P$  which changes places of $a$ and $c$. So, we obtain the third local symmetric space as follows
$$\sll(2k)/\gl{k}.$$

{\bf Conclusion.} We obtain the $\z$-symmetric spaces corresponding to the pairs
$(\g,\g_e)$ given in the following table:

\begin{center}
\begin{tabular}{|c|c|}
\hline
&\\
$\g$ & $\g_e$ \\
&\\ \hline\hline
&\\
$\sll(2k)$ & $\gll(k)$\\
& \\
\hline
&\\
$\sll(n)$ & $\soo(n)$\\
& \\
\hline
&\\
$\sll(2m)$ & $\spp(2m)$\\
& \\
\hline
&\\
$\sll(k_1+k_2)$ & $\soo(k_1)\oplus\soo(k_2)$\\
& \\
\hline
&\\ 
$\sll(2(k_1+k_2))$ & $\spp(2k_1)\oplus\spp(2k_2)$\\ 
& \\ \hline
\end{tabular}

\bigskip

\textit{Table 4}
\end{center}

\medskip

We summarize the result obtained in Section \ref{scl} as follows.

\begin{theorem}\label{tmain} All local complex $\z$-symmetric spaces in cases $A_l,\: l\ge 1$, $B_l,\: l\ge 2$, $C_l,\: l\ge 3$ or $D_l,\: l\ge 4$ are given in Tables 1, 2, 3, and 4. Each space $\g/\g_e$ is uniquely, up to a weak equivalence of respective $\z$-grading, defined by $\g$ and $\g_e$.
\end{theorem}

\begin{center}
\textbf{Acknowledgement}
\end{center}

The authors are grateful to the referee for a number of valuable comments and a correction in one case.

\end{document}